# Linéarisation d'une itération inversible bornée dans $\mathbb{R}^d$ par des fonctions de Weierstrass

**G. Cirier**
LSTA. Université Paris VI, France
Email: guy.cirier@gmail.com

**Résumé**

Soit une itération de $\mathbb{R}^d$ dans $\mathbb{R}^d$ définie par un difféomorphisme polynomial borné. On montre que les courbes semi invariantes tendent vers des courbes paramétrées par des fonctions de Weierstrass. Cela justifie les calculs d'échelle d'autosimilarité et de dimension fractale comme le pratiquent des praticiens sur des itérations chaotiques. On applique ces résultats au calcul différentiel.

**Abstract**

In this paper, we study an iteration in $\mathbb{R}^d$ defined by a diffeomorphism polynomial bounded. Semi invariant curves tend to curves with parametric Weierstrass-Mandelbrot's functions. So, self-similarity and fractal dimension are justified. We apply these results to partial differential calculus.

**Mots clefs** : Itération bornée, courbes semi invariantes, fonctions de Weierstrass, Hénon, EDP, Navier-Stokes.
**Key words**: Bounded iterations, semi invariant curves, Weierstrass functions , Hénon, PDE, Navier-Stokes.

## I - Introduction

Soit $f$ une application de $\mathbb{R}^d$ dans $\mathbb{R}^d$. On appelle itération la même application $f$ si elle est itérée indéfiniment : $f \circ f \circ ... \circ f = f^{(p)}$. L'équation de Perron Frobenius qui caractérise les mesures invariantes sous l'itération $f$, a été le point de départ de deux articles précédents : le comportement déterministe de $f$ prévaut dans $[5]$; mais, dans $[6]$, le point de vue probabiliste diffère. On approfondit ici l'approche déterministe en la situant dans un ensemble borné $C \subset \mathbb{R}^d$ tel que $f(C) \subset C$. Les résultats obtenus ici justifient certaines méthodes utilisées par les praticiens.

### 1- Hypothèse H0

*Les itérations $f$ et $f^{-1}$ sont polynomiales bornées dans* C *ayant un point fixe $f(0)=0$, toutes les v.p. $\lambda > 1$ sont réelles positives non résonantes : $1 \neq \lambda^n$ pour tout $n \in \mathbb{Z}^d$.*
Si les v.p. réelles $|\lambda| > 1$, sont négatives, quitte à itérer 2 fois $f$, $\boldsymbol{a}(\lambda^2 \boldsymbol{t}) = f \circ f(\boldsymbol{a}(\boldsymbol{t}))$, **on travaillera** ***donc toujours avec*** $\lambda > 1 \in \mathbb{R}^{+q}$.

### 1- Résultats connus et admis : si $f$ est un $C^k$ - difféomorphisme dans $\mathbb{R}^d$ , $k \geq 1$

Sous ces conditions, même si $f$ n'est qu'un $C^k$ - difféomorphisme avec $k \geq 1$, Steinberg$[10]$a construit, $d$ ***fonctions de linéarisation*** $C^k$ *:* $\varphi_\ell \circ f(\boldsymbol{a}) = \lambda_\ell \varphi_\ell(\boldsymbol{a})$, $\ell = 1,..,d$ au voisinage de $0$.
De même, on peut établir une bijection entre l'itération $\boldsymbol{a}_1 = f(\boldsymbol{a}) \subset C$ à partir d'un point $\boldsymbol{a}_0 \in C$ et l'équation fonctionnelle $\boldsymbol{a}(\lambda \boldsymbol{t}) = f(\boldsymbol{a}(\boldsymbol{t})), \boldsymbol{t} \in \mathbb{R}^d$. On obtient $d$ ***fonctions semi invariantes*** $\boldsymbol{a}(\boldsymbol{t})$: $a_\ell(\lambda \boldsymbol{t}) = f_\ell(\boldsymbol{a}(\boldsymbol{t})$ $\ell = 1,..,d$ . Dans les deux cas, les applications $\varphi$ et $\boldsymbol{a}$ sont construites localement avec des séries. En outre, on a l'unicité de $\varphi$ et de $\boldsymbol{a}$ pour une matrice jacobienne unité en $0$.

On construit d'abord $a(t)$ avec les séries de Steinberg définies au voisinage de 0. Soit rayon de convergence $\rho > 0$ de $a(t)$. Puis, en choisissant $t : |t| \leq \rho_0 < \rho$, par itérations successives, on prolonge $a(\lambda^p t) \in C$ pour $\forall p \in Z$ , donc pour $\forall t \in \mathbb{R}^d$. Cette solution unique sera très utile ici.

Comme les coordonnées de $t$ de $f^{(p)}$ relatives aux $|\lambda| < 1$ tendent vers 0 quand $p \to \infty$, on ne s'intéresse ici qu'aux $q$ v.p. $\lambda \in \mathbb{R}^{+q} \| \lambda | > 1$ . On leur associe $t \in \mathbb{R}^q$.

**Soit E l'équation $a(\lambda t) = f(a(t)) | \lambda > 1, t \in \mathbb{R}^q$. On cherche ici les solutions de E. Par exemple, une itération qui commence avec $a(t_0, t'_0)$ où $0 \neq t'_0 \in \mathbb{R}^{d-q}$, ne vérifiera pas en général E. Mais, par continuité, $a_0 = a(t_0, 0) = \lim_{p \to \infty} a(t_0, \lambda'^p t'_0)$ vérifie E, sans dépendre de $t'_0$, donc de $a(t_0, t'_0)$.**

Enfin, si l'unicité du point fixe $0$ dans $C$ n'est pas acquise, on peut encore rencontrer des difficultés qui ne sont pas étudiées ici. Les solutions trouvées sont locales.

## II- Fonctions de Bohr presque périodiques dans un Banach [1] et application aux itérations

Soit P.P. l'ensemble des fonctions presque periodiques bornées dans $\mathbb{R}^d$. Soit $a(t), t \in \mathbb{R}^q, q \leq d$ .

### 1- Définition de Bochner (rappel [4],[3])

*Une fonction $a(t) \in P.P.$ au sens de Bohr si et seulement si l'ensemble des fonctions translatées par tout $c$ $\{a(t+c), t \in \mathbb{R}^q\}$ est relativement compact pour la topologie de la convergence uniforme.*

On dispose alors d'un arsenal de résultats proches de ceux des séries de Fourier [1]. Le théorème d'approximation est démontré dans des espaces de Banach. Il est utilisé ici dans $\mathbb{R}^d$.

### 2- Théorème 1 (rappel de l'approximation de Bohr [1])

*Soit un polynôme trigonométrique : $P_n(t) = \Sigma_{k=1}^n c_k e^{i\mu_k t} \in \mathbb{R}^d$ où $(c_k \in \mathbb{R}^d, \mu_k \in \mathbb{R}^q)$ et $\mu_k t$ le produit scalaire. Si $a(t) \in P.P.$, il existe une suite de $P_n(t)$ tels que pour tout $\varepsilon > 0$ : $\sup_{t \in \mathbb{R}^q} \|a(t) - P_n(t)\| < \varepsilon$.*

*Si $c(\mu) = \lim_{T \to \infty} \frac{1}{(2T)^q} \int .. \int_{-T}^{T} a(t) e^{-i\mu t} \, dt$ est vecteur coefficient de Fourier-Bohr, alors $\Lambda(a) = \{\mu \in \mathbb{R}, c(\mu) \neq 0\}$ est au plus dénombrable. $a(t)$ est uniformément continue et développable en série de Fourier : $a(t) \sim \sum_{\mu \in \Lambda} c(\mu) e^{i\mu t}$. Cette représentation est unique.*

*La convergence a lieu en moyenne quadratique et l'on a l'égalité de Perceval.*

On remplace ici le groupe additif utilisé par Bohr par le groupe multiplicatif dans $\mathbb{R}^q$.

### 3-Théorème 2

*Sous H0, $\{a(t), t \in \mathbb{R}^q\}$ et $a(ct)$ pour tout $c$ sont relativement compacts et $a(t) \in P.P.$*

Soit $A_N = \{(a(\lambda^p t); \|t\| \leq \rho_0, p \geq N\}$. Son adhérence $\widehat{A}$ est fermée et contenue dans C borné, donc compacte. Donc l'image de $a(t)$ est relativement compacte. « Translatons » $t$ par $c$ : $ct$ de sorte que $a(ct) \in C$. Soit $n$ assez grand pour que $c = \alpha \lambda^n$ avec $\|\alpha\| < 1$ et $\|\alpha t\| < \|t\| < \rho_0$, $ct = \alpha t \lambda^n$. Comme $f$ est inversible, $f^{-(n)} a(ct) = a(\alpha t)$ avec $|\alpha t| \leq \rho_0$ qui redonne $A_N$ d'adhérence $\widehat{A}$. Donc $a(t) \in P.P.$.

### Corollaire 1

*Sous H0, les coordonnées du graphe de $a(t)$ dans $\mathbb{R}^d$ sont approximées uniformément par un paramétrage presque périodique en $t \in \mathbb{R}^q$.*

**II- Solution de l'équation E : $a(\lambda t) = f(a(t)), t \in \mathbb{R}^q$ dans $P.P.$**

Soit la transformée de Fourier de $f(a(t))$ pour $\mu$ : $F_\mu(f(a(t)) = \lim_{T\to\infty} \frac{1}{(2T)^q} \int ... \int_{-T}^{T} f(a(t)e^{-it\mu} dt$

La formule de Taylor en dimension $d$ pour $f$ en somme de monômes est: $f(a) = \Sigma_{n=1}^{n=q} f^n(0) a^n / n!$.

On note $\lambda$ ces valeurs propres et $|\lambda| = \Pi_{j=1}^{j=q} \lambda_j$ leur produit.

### 1- Lemme de linéarisation

*Sous H0, les presque périodes (p.p.) $\mu \in \Lambda$ solutions possibles de l'équation* **E** *sont 0 et $\mu = \lambda^k |\ k \in \mathbb{Z}$. On écrira $a(t) = \Sigma_\mu c(\mu) e^{i\mu t} + c(0)$*

Soient 0 (si $c(0) \neq 0$) et $\mu > 0$ des p.p. possibles. Si $a_i(t)$ est une coordonnée quelconque de $a(t)$, il suffit d'examiner un des monômes $\Pi_i (a_i(t)^{n_i}$ de $a(t)^n$ avec $\Sigma n_i = n \leq q$. Comme on a : $a_i(t) \sim \Sigma_\mu c_i(\mu) e^{i\mu t} + c_i(0)$ avec $\mu > 0 \in \Lambda$, le monôme devient : $\Pi_i (a_i(t)^{n_i} = \Pi_i (\Sigma_\mu c_i(\mu) e^{i\mu t} + c_i(0))^{n_i}$.
Avec la formule de la multinomiale, on écrit alors:

$\Pi_i (a_i(t))^{n_i} = \Pi_i (n_i! \Sigma_j \Sigma_\mu \Pi_{m_j=1}^{n_i} (c_i(0)^{s_j} (c_i(\mu)^{m_j} / m_j! s_j!) e^{it\Sigma_j m_j \mu}$. On applique alors Fourier $F_{\mu_0}$ :
$F_{\mu_0}(a(\lambda t)) = c(\mu_0 / \lambda) / |\lambda| = F_{\mu_0}(f(a(t)))$.

Si $\mu_0 = 0$, on devra avoir pour le monôme de $f$ : $\Sigma_{j\mu} m_j \mu = 0$, donc $\mu = 0$.

Si $\mu_0 > 0$, on compare alors les exposants de $e^{i\mu t}$ à $\mu_0$. On obtient 2 relations : $\mu / \lambda = \mu_0$, et $\Sigma_{j\mu} m_j \mu = \mu_0$. La première implique, puisque $f$ est inversible, $\mu = \mu_0 \lambda^k |\ k \in \mathbb{Z}$. On a une relation d'équivalence entre les $\mu > 0$ modulo les puissances de $\lambda$. Sous H0, la seconde relation devra être vraie pour des $m_j \geq 2$ avec $f$ et $f^{-1}$. Donc, pour toute itération et tout $m \geq 2$, $\mu_m \leq \mu_0 / 2$. En résumé, les seules solutions positives possibles seront de la forme $\mu_m \lambda^k |\ k \in \mathbb{Z}, \mu_m \leq \mu_0 / 2$.

Montrons l'unicité de $\mu_0$ modulo $\lambda$. En effet, si $\mu_m \leq \mu_0 / 2$ et $\mu_m \neq \mu_0 \lambda^\ell$ est une autre solution, on aura : $\Sigma_j m_j \mu_m \lambda^k + \Sigma_j m"_j \mu_0 \lambda^{k"} = \mu_m \lambda^{k'}$. Ce qui donne avec $r = \mu_0 / \mu_m \geq 2$, $\Sigma_j m_j \lambda^{k-k'} + r \Sigma_j m"_j \lambda^{k"-k'} = 1$. Ce qui est impossible. On prend $\mu_0 = 1$ et l'on pose $\mu = \lambda^k |\ k \in \mathbb{Z}$.
Ceci est vrai pour toutes les coordonnées de $\lambda$. Précisons les coefficients relatifs à ces p.p.

### 2 - Théorème de linéarisation 3

*Sous H0, les coefficients de $a(t) \in P.P.$ vérifiant* **E** *sont $a(t) = \Sigma_\mu c(\mu) e^{i\mu t} + c(0)$ où :*

*- $c(0)$ est point fixe non nul de $|\lambda| f$ dans $C \subset \mathbb{R}^d$ (condition essentielle pour que 0 soit une p.p.) ;*

*- $c(\mu)$ vérifie la récurrence linéaire $c(\mu / \lambda) / |\lambda| = \frac{\partial f(c(0))}{\partial a} c(\mu)$. En posant $c(\mu) = cr^k$, on a :*

$$\frac{\partial f(c(0))}{\partial a} c = \rho c$$

*$c_1$ est donc vecteur propre pour la valeur propre $\rho = 1 / (r|\lambda|)$.*

*$a(t) = c_1 \Sigma_{i=1}^{i-d} w_i(t) + c(0) + c_1 \Sigma_{i=1}^{i-d} \Sigma_k r_i^{|k|}$ est somme des fonctions de Weierstrass $w_i(t) = \Sigma_k r_i^{|k|} (1 - e^{i\lambda^k t})$, chacune dépendant de $r_i$ défini à partir d'une v.p. $\rho_i$.*

On va obtenir d'abord les coefficients de Fourier pour $\mu = 0$ : $c(0) / |\lambda| = f(c(0)$ . $c(0)$ est point fixe de $|\lambda| f$ dans $C \subset \mathbb{R}^d$. Pour que $0$ soit une p.p., il faut que $c(0) \neq 0$.

La transformée de Fourier $\boldsymbol{F}_\mu$ de $\boldsymbol{a}(\lambda \boldsymbol{t}) = f(\boldsymbol{a}(\boldsymbol{t}))$ pour une période $\mu = \lambda^k \mid k \in \mathbb{Z}$ donnera un résultat nul pour tout multiple entier $m \geq 2$ de $\mu$. Donc, le résultat doit être linéaire en $\boldsymbol{c}(\mu)$: $\boldsymbol{F}_\mu f(\boldsymbol{a}(\boldsymbol{t})) = A(f)\boldsymbol{c}(\mu)$ où $A(f)$ est une matrice constante. Précisons $A(f)$ avec le calcul du monôme de la formule de Taylor transformé par Fourier:
$(\boldsymbol{a})^{\alpha_n} = a^{\alpha_1}...a^{\alpha_i} // \alpha_1!...\alpha_i!$ donne $\boldsymbol{F}_\mu(\boldsymbol{a}(\boldsymbol{t}))^{\alpha_n} = c_1(\boldsymbol{0})^{\alpha_1} / \alpha_1!...c_{i-1}(\boldsymbol{0})^{\alpha_{i-1}} / \alpha_{i-1}! \boldsymbol{F}_\mu a_i(t)^{\alpha_i} / \alpha_i!$
Or, $a_i(\boldsymbol{t}) = \Sigma_\mu c_i(\mu) e^{i\mu t} + c_i(0)$ , donc : $\boldsymbol{F}_\mu a_i(t)^{\alpha_i} / \alpha_i! = \alpha_i (c_i(0)^{\alpha_i - 1} / \alpha_i!) c_i(\mu)$ qui n'est autre que la dérivée du monôme par rapport à $c_i(0)$. Par sommation, on obtient : $\boldsymbol{F}_\mu f(\boldsymbol{a}(\boldsymbol{t})) = \frac{\partial f(\boldsymbol{c}(\boldsymbol{0}))}{\partial \boldsymbol{a}} \boldsymbol{c}(\mu)$.

**Corollaire 2**

*Lorsque l'on itère indéfiniment la fonction $f$ , $\boldsymbol{a}(\boldsymbol{t})$ se comporte asymptotiquement comme la fonction de Weierstrass relative à la plus petite v.p. $r_0$ . Elle devient donc auto similaire.* Si $|r_0| = \min_i |r_i|$ , alors

$$\boldsymbol{a}(\boldsymbol{t}) \to \boldsymbol{c} w_{r_0}(t)$$

*Si $\lambda$ est unidimensionnel,* l*a dimension de $\boldsymbol{a}(\boldsymbol{t})$ est asymptotiquement (par la méthode des boites)* :

$$D = 2 + \log|r_0| / \log|\lambda|$$

Itérons $p$ fois $f$ . Tout se passe comme si $\boldsymbol{a}(\boldsymbol{t})$ devenait $\boldsymbol{a}(\lambda^p \boldsymbol{t})$ . Comme $\boldsymbol{a}(\boldsymbol{t})$ est somme de fonctions de Weierstrass, chacune dépendant d'une racine, on observe une invariance d'échelle distincte pour chacune: $w_i(\lambda^p t) = r_i^{-|p|} w_i(t)$. Si $|r_0| = \min_i |r_i|$, mettons en facteur $r_0^{-|p|}$. Alors $(r_0 / r_i)^{|p|} w_i(t) \to 0$ uniformément quand $p \to \infty$ . Il en est de même pour les termes constants.
On connaît ensuite la dimension fractale de $w_r(t)$ en identifiant $|r_0|^k = |\lambda|^{(D-2)k}$ , donc celle de $\boldsymbol{a}(t)$ .

Si $r_0$ est complexe, $r_0 = \rho e^{i\theta}$ ; $w_r(t) = 2\Sigma_k \rho^k \cos k\theta (e^{i\lambda^k t} - 1)$

**Conditions de convergence de Hardy [9]**

*Rappelons que la série de Weierstrass* $w_r(t) = \Sigma_{k=-\infty}^{\infty} (1 - e^{i\lambda^k t}) r^{|k|}$ *ne converge que si* $|\lambda r| > 1 \,\&\, |r| < 1$ . Ici, *pour converger, la v.p.* $\rho = 1/(r|\lambda|)$ *doit vérifier* $|\lambda| \rho > 1 \,\&\, \rho < 1$ .

**3- Remarques**

- On peut remplacer les fonctions polynomiales par des fonctions analytiques ( Favard [7]) sans changer les résultas: si l'application est un difféomorphisme $C^\infty$, à l'intérieur d'un disque de convergence $\|\boldsymbol{a}\| < \rho_1$, $f$ est approximée uniformément par un polynôme. Par suite, le rayon de convergence de la série $\boldsymbol{a}(\boldsymbol{t})$ de Steinberg étant $\rho$ , il suffit de prendre $\boldsymbol{t}$ suffisamment petit pour que $\|\boldsymbol{a}(\boldsymbol{t})\| < \rho_1$ et $\|\boldsymbol{t}\| < \rho$ pour conserver l'ensemble des résultats.
- Si les v.p. sont complexes, $\lambda = \rho e^{i\theta}$ . On peut approximer $\theta$ avec Dirichlet par $\theta = 2k\pi / p$ . Si $\theta = 2k\pi / p$ , en itérant $f^{(p)}$, $\lambda^p = \rho^p > 1$ , on retrouve les conditions de la démonstration. Mais, ce faisant, on introduit les points de cycles d'ordre $p$ . La solution trouvée en un point du cycle doit être itérée $k = 1,...,p$ fois. C'est ce qui peut se produire si $\lambda < -1$ . Ces solutions ont une fréquence $1/p$ .

**4- Conséquences**

On a une justification mathématique des constats et des méthodes utilisées par les praticiens :
- Comme la fonction p.p. tend asymptotiquement vers une fonction de Weierstrass, le constat de l'auto similarité est justifié, ainsi que les méthodes de dimension fractales. L'auto corrélation a aussi été utilisée pour étudier la fonction de Weierstrass et s'applique aussi.
- La sensibilité aux conditions initiales ne joue pas sur la forme de la courbe. Il semble les chercheurs aient été perturbés par la « périodicité multiplicative » et ont attribué un peu rapidement ce désordre apparent aux conditions initiales. Le « vrai » désordre est dû aux fonctions de Weierstrass.

- En fait, on peut approximer comme l'on veut la forme générale de la courbe par des techniques mathématiques usuelles sans pouvoir situer exactement un point à un moment donné sur sa trajectoire. Ainsi, la dimension fractale est liée à l'invariance d'échelle $r_0$. Que se passe-t-il si $\lambda$ n'est pas unique ? Quelle est alors la dimension fractale du paramétrage?

**III – Exemple : cas de Hénon borné**

L'itération est définie par $(a_1,b_1)=\gamma a+b+h(a),\beta a)$ où $\gamma=\lambda+\lambda'$ et $\beta=-\lambda\lambda'$ avec $|\lambda|>1$ et $|\lambda'|<1$. Sa courbe semi invariante se paramètre avec $\boldsymbol{a}(t)=(a(\lambda t),a(t))$ vérifiant: $a(\lambda t)=\gamma a(t)+\beta a(t/\lambda)+h(a(t))$ . Si $\lambda$ est négatif, en 2 itérations, on se retrouve dans le cas positif. Mais, la récurrence à étudier est alors celle de l'itération répétée deux fois et l'on a des cycles d'ordre 2 au lieu de points fixes. On a :

$$(a_2,b_2)=\gamma(\gamma a+b+h(a))+\beta(\gamma a+b+h(a))+h(\gamma a+b+h(a)),\beta(\gamma a+b+h(a))$$

La v.p. plus grande que 1 est $\lambda^2$.

**Proposition**

*Sous H0, et sous les conditions de convergence de Hardy* **H** *énoncées ci-dessous, le paramétrage* $a(t)$ *de la courbe d'Hénon bornée* $\boldsymbol{a}(t)$ *est la somme de deux fonctions de Weierstrass de presque période* $\mu_k=\lambda^{2k}$ *et de coefficients* $r$ *et* $r'$ *si* $\lambda<-1$. $w_r(t)$ et $w_{r'}(t)$ *sont approximées par des polynômes trigonométriques :*

$$P_n(t)=\Sigma_{k=-n}^{n}(1-e^{i\lambda^{2k}t})(\alpha r^{|k|}+\alpha' r'^{|k|})$$

Calculons $c(0)=(c_1,c_2)$. Si l'on a un cycle d'ordre 2, on aura pour l'autre point du ycles lesmêmes équations, à une translation près.

$$\begin{pmatrix} c_1/\lambda^2 \\ c_2/\lambda^2 \end{pmatrix}=\begin{pmatrix} \gamma(\gamma c_1+c_2+h(c_1))+\beta(\gamma c_1+c_2+h(c_1))+h(\gamma c_1+c_2+h(c_1)) \\ \beta(\gamma c_1+c_2+h(c_1)) \end{pmatrix}$$

Calculons ensuite $c(\mu)=(d_1r^k,d_2r^k)$ avec les équations aux valeurs propres :

$$\begin{pmatrix} \gamma+h'(\gamma c_1+c_2+h(c_1)) & 1 \\ \beta & 0 \end{pmatrix}\begin{pmatrix} \gamma+h'(c_1) & 1 \\ \beta & 0 \end{pmatrix}\begin{pmatrix} d_1 \\ d_2 \end{pmatrix}=\rho\begin{pmatrix} d_1 \\ d_2 \end{pmatrix}$$

Le déterminant doit être nul avec $\rho=1/\lambda^2 r$ :

$$\begin{vmatrix} (\gamma+h'(\gamma c_1+c_2+h(c_1)))(\gamma+h'(c_1)-\rho & \gamma+h'(\gamma c_1+c_2+h(c_1)) \\ \beta(\gamma+h'(c_1)) & \beta-\rho \end{vmatrix}=0$$

qui fournit deux racines réelles $r$ et $r'$ si :

$$((\gamma+h'(\gamma c_1+c_2+h(c_1)))(\gamma+h'(c_1)+\beta)^2+4(\gamma+h'(\gamma c_1+c_2+h(c_1))\beta(\gamma+h'(c_1))\geq 0$$

Sinon, on peut avoir des racines complexes. Notons que si on peut linéariser la récurrence avec $\lambda>1$, on a une relation bien plus simple : $c(\mu_k/\lambda)/\lambda=(\gamma+c_0)c(\mu_k)+\beta c(\mu_k\lambda)\lambda$ qui induit une linéarité cyclique.

**Nota** : Dans de l'itération classique d'Hénon, on a $h(a)=-\sigma a^2$ avec $\sigma=1{,}4$ et $\beta=0{,}3$ avec $\gamma=-2\sigma*0{,}6313=-1{,}7678$. Les v.p. réelles au point fixe 0 sont $\lambda\approx-1{,}9237$ et $\lambda'\approx 0{,}1559$. La courbe $(a,b)=(a(t\lambda),a(t))$ définie par le paramétrage $a(\lambda^2 t)=a(t)/r=w_r(t)/r$ connu est: soit $(1-r\beta)a(t\lambda)=\gamma a(t)+h(a(t))$, soit $(1-r\beta)a(t)/r=\gamma a(t\lambda)+h(a(t\lambda))$. Cette dernière formule génère une infinité de branches dont on peut en contrôler les approximations.

**IV- Application aux équations différentielles ordinaires ou aux dérivées partielles**

On considère les équations de la forme $\partial\boldsymbol{a}/\partial\boldsymbol{x}=F(\boldsymbol{a})$. Les inconnues forment un vecteur $\boldsymbol{a}\in\mathbb{R}^d$. La variable $\boldsymbol{x}=(\boldsymbol{q},t)\in\mathbb{R}^k$ s'interprète comme la position $\boldsymbol{q}=(q_1,...q_{k-1})\in\mathbb{R}^{+k-1}$ supposée positive de $\mathbb{R}^{+k-1}$ et le temps $t_k = t \in \mathbb{R}^+$. $F(\boldsymbol{a})$ est une matrice $(k,d)$ de polynômes en $\boldsymbol{a} \in \mathbb{R}^d$ dans $\mathbb{R}^d$.
Si $1 < k < d$, plusieurs coordonnées de $\boldsymbol{a}$ sont fonction des mêmes variables. $k = 1$ est le cas d'une EDO avec pour seule variable $t$. La référence sur la question reste Arnold$[2]$. On privilégie ici $k = d$.
Pour étudier cette équation différentielle, il faut d'abord la traduire en termes d'itération.

**1-Définition et notations**

*On prend $k = d$. On appelle itération différentielle l'application $f(\boldsymbol{a},\delta)$ de $\mathbb{R}^d$ dans $\mathbb{R}^d$, définie par $\boldsymbol{a}_1 = f(\boldsymbol{a},\delta) = \boldsymbol{a} + \delta F(\boldsymbol{a})$, $\delta \in \Delta = \mathbb{R}^{+d} \cap \{\delta_0 \geq \delta \geq 0\}$.*

On relie $\delta$ à $\boldsymbol{x}$ par la relation $\delta = \boldsymbol{x}/n$, $n \in \mathbb{N}, \boldsymbol{x} \in \mathbb{R}^d$. Les zéros réels de $F$ sont les points fixes de $f(\boldsymbol{a},\delta)$ et ne dépendent pas de $\delta \in \Delta$. On posera $\boldsymbol{x}F = x_i F_i; i = 1,..,d.$ On écrit l'itération dans la base des vecteurs propres de $\partial F(\boldsymbol{0})/\partial \boldsymbol{a}$, de sorte que si $\rho_i$ est valeur propre de $\partial F(\boldsymbol{0})/\partial \boldsymbol{a}$, on a : $\lambda_i = 1 + \delta_i \rho_i = 1 + x_i \rho_i / n$. On associe ainsi de façon unique à $\rho_i$, les quantités $\lambda_i$, $x_i$, $a_i$ et $a_{1i}$.

**Hypothèse1**

*Sous H0, pour tout $\delta = \boldsymbol{x}/n \leq \delta_0$, $f(\boldsymbol{a},\delta)$ applique un ensemble borné $C$ de $\mathbb{R}^d$ dans lui-même. On a 1 v.p. $\rho > 0$ qui rend $\lambda > 1$, tous les autres $\rho_i$ sont négatifs.*

On associe à cette v.p. $\rho > 0$ la variable $t$ et, aux autres $\rho_i < 1$, le vecteur de position $\boldsymbol{q}$. On choisit $n$ assez grand pour que $\sup|\delta_0 \rho_i| < \varepsilon$.

Au voisinage d'un point fixe $\boldsymbol{0}$, pour tout $\delta \leq \delta_0$ suffisamment petit, $f(\boldsymbol{a},\delta)$ est inversible, donc un $C^\infty$ difféomorphisme. D'après le théorème 4, on peut faire les calculs avec une seule v.p. positive puisque seule la partie périodique est modifiée quand on a plusieurs v.p. plus grandes que 1et l'on remplace $\lambda$ par $|\lambda|$.

**2- Proposition 1**

*Sous H1, l'itération différentielle admet une solution asymptotique presque périodique quand $n \to \infty$. Pour $n$ fixé, $\lambda = 1 + t\rho/n$ est la seule v.p ; $\lambda > 1$ $t$. La solution s'écrit :*

$$\boldsymbol{a}(u) \sim \sum_{\mu \in \Lambda} \boldsymbol{c}(\mu) e^{i\mu u} + \boldsymbol{c}(0) \text{ avec } \mu = \lambda^k \text{ et } 0.$$

*Les $\boldsymbol{c}(0)$ vérifient $\boldsymbol{c}(0)/\lambda = \boldsymbol{c}(0) + \delta F(\boldsymbol{c}(0))$ et les $\boldsymbol{c}(\mu)$ vérifient la relation linéaire :*

$$\boldsymbol{c}(\mu/\lambda)/\lambda = \boldsymbol{c}(\mu) + (\boldsymbol{x}/n)\partial F(\boldsymbol{c}(0))/\partial a\, \boldsymbol{c}(\mu)).$$

*En prenant $\boldsymbol{c}(\mu_k) = (r)^k \boldsymbol{c}$* avec $1/r = (1 + \sigma t/n)\lambda$, l*a solution est :*

$$\sigma t\, \boldsymbol{c} = \boldsymbol{x} \partial F(\boldsymbol{c}(0))/\partial a\, \boldsymbol{c}$$

*$\boldsymbol{c}$ est vecteurs propre de $\boldsymbol{x}\partial F(\boldsymbol{c}(0))/\partial a$ pour la v.p. $\sigma t$.*

Si $n \to \infty$, $\boldsymbol{c}(0) = -\boldsymbol{x}F(\boldsymbol{c}(0))/\rho t$ est point fixe de $-\boldsymbol{x}F/\rho t$ et :

$$\boldsymbol{a}(u) \to ((\sigma + \rho)t)^{-1}(1 - e^{-t(\sigma+\rho)}) w(u) \boldsymbol{c}$$

*$w(u) = \Sigma_k e^{-t(\sigma+\rho)|k|}(1 - \exp(ie^{t\rho k}u))$ est la fonction de Weierstrass correspondant au plus petit $e^{-t(\sigma+\rho)}$. La condition de convergence de Hardy est alors : $\sigma < 0$.*
*La dimension est alors $D = 1 - \sigma/\rho$*

Ici encore, on construit au voisinage de $0 \in C$ l'équation d'itération: $\boldsymbol{a}(\lambda \boldsymbol{u}) = f(\boldsymbol{a}(\boldsymbol{u}),\delta) = \boldsymbol{a}(\boldsymbol{u}) + \delta F(\boldsymbol{a}(\boldsymbol{u}))$. On a une seule v.p. plus grande que 1 : $\lambda = 1 + \rho t/n > 1$, toutes les autres $\lambda_i < 1$. Quand $p \to \infty$, toutes les coordonnées de $\boldsymbol{u}$ tendent vers 0 sauf celle relative à $t$ que l'on notera encore $u$. On a alors une représentation de $f(\boldsymbol{a}(u),\delta)$ par des fonctions p.p. : $\boldsymbol{a}(u)$ est développable en série de Fourier : $\boldsymbol{a}(u) \sim \sum_{\mu \in \Lambda} \boldsymbol{c}(\mu) e^{i\mu u} + \boldsymbol{c}(0)$, $\Lambda(a) = \{\mu \in \mathbb{R}, c(\mu) \neq 0\}$ est au plus dénombrable.

Calculons d'abord $\boldsymbol{c}(0)$ : $\boldsymbol{c}(0)/\lambda = \boldsymbol{c}(0) + \delta F(\boldsymbol{c}(0))$. Ce qui donne, quand $n \to \infty$ : $-\rho t\boldsymbol{c}(0) = \boldsymbol{x}F(\boldsymbol{c}(0))$. $\boldsymbol{c}(0)$ est point fixe non nul de $-\boldsymbol{x}F/\rho t$. Calculons les coefficients $\boldsymbol{c}(\mu)$ avec le théorème 3 en posant $\lambda_i = 1 + \rho_i x_i / n$

$\boldsymbol{c}(\mu/\lambda)/\lambda = \boldsymbol{c}(\mu) + (\boldsymbol{x}/n)\partial F(\boldsymbol{c}(0))/\partial a\boldsymbol{c}(\mu))$

La récurrence en $\boldsymbol{c}(\mu)$ devient linéaire avec $\mu = \lambda^k$. En posant : $\boldsymbol{c}(\mu) = (r)^{|k|}\boldsymbol{c}$ et $1/\lambda r - 1 = \sigma t/n$, c'est-à-dire $1/r = (1 + \sigma t/n)\lambda$, on a l'équation linéaire ;

$\sigma t\, \boldsymbol{c} = \boldsymbol{x}\partial F(\boldsymbol{c}(0))/\partial a\, \boldsymbol{c}$ où $\boldsymbol{c}$ est vecteur propre de $\boldsymbol{x}\partial F(\boldsymbol{c}(0))/\partial a$ pour la v.p. $\sigma t$. $\boldsymbol{a}(u)$ s'écrit alors: $\boldsymbol{a}(u) = \alpha \Sigma_k \boldsymbol{c} r^k \exp(i\lambda^k u) + \boldsymbol{c}(0)$

Posons $|k| = |k'|n + s, s = 1,..,n-1, k' \in \mathbb{Z}$, avec $\alpha = 1/n$ et $v = s/n$. Si $n \to \infty$ :

$r^{|k|} \to \exp - t(\rho + \sigma)(|k'| + v)$ et $\lambda^k \to \exp(\rho t(|k'| + v))$

$$\Sigma_k \boldsymbol{c}(\mu)e^{i\mu u} = \Sigma_r (1/n)\boldsymbol{c}\Sigma_{s=1}^{n-1} e^{-t(\sigma+\rho)(|k'|+v)} \exp(ie^{t\rho(k'+v)}u) \to \Sigma_q \boldsymbol{c}\int_0^1 e^{-t(\sigma+\rho)(|k'|+v)} \exp(ie^{\rho(k'+v)}u)dv$$

Itérons la relation différentielle $np_0$ fois, $u$ devient $(1 + \rho t/n)^{np_0} u \to ue^{p_0\rho t}$, donc, $\exp(ie^{t\rho(k'+v)}u)$ devient $\exp(ie^{t\rho(k'+p_0+x)}u)$ où $p_0$ est arbitrairement grand et $v$ devient négligeable devant $p_0$ ( $p_0 + v \in (p_0, p_0 + 1)$ ) :

$$\int_0^1 e^{-t(\sigma+\rho)|k'|+v)} \exp(ie^{\rho(k'+v+p_0)}u)dv \to ((\sigma + \rho)t)^{-1}(1 - e^{-t(\sigma+\rho)})e^{-t(\sigma+\rho)k'} \exp(ie^{\rho(k'+p_0)}u).$$

On peut alors identifier $k' + p_0$ à $k$.

**3- Solution asymptotique quand $t \to \infty$**

Il est naturel de chercher la solution asymptotique de ces équations quand $t \to \infty$ alors que les variables de position $\boldsymbol{q}$ restent bornées. Comme on fait jouer un rôle dissymétrique à $t$ et à $\boldsymbol{q}$, on note alors $\boldsymbol{a}_{1t} = \boldsymbol{a}_t + (t/n)F_t(\boldsymbol{a})$ et $\boldsymbol{a}_{1q} = \boldsymbol{a}_t + (\boldsymbol{q}/n)F_q(\boldsymbol{a})$. Les coefficients de Fourier-Bohr s'écrivent :

$c_t(0)/\lambda = c_t(0) + (t/n)F_t(\boldsymbol{c}(0))$ et $\boldsymbol{c}_q(0)/\lambda = \boldsymbol{c}_q(0) + (\boldsymbol{q}/n)F_q(\boldsymbol{c}(0))$. Si $\boldsymbol{q}/t \to 0$, on obtient :

$-\rho c_t(0) = F_t(\boldsymbol{c}_t(0), 0_q)$ et $\boldsymbol{c}_q(0) = 0$.

Calculons les coefficients $\boldsymbol{c}(\mu)$ avec $\mu = \lambda^k$ qui rend la récurrence en $\boldsymbol{c}(\mu)$ linéaire. En posant $\boldsymbol{c}(\mu) = (r)^k \boldsymbol{c}$, on a :

$-\rho t\, \boldsymbol{c}/\boldsymbol{r} = \boldsymbol{x}\partial F(\boldsymbol{c}(0))/\partial a\, \boldsymbol{c}$. Si $\boldsymbol{q}/t \to 0$, le déterminant de la relation devient : $(-\rho/\boldsymbol{r}\lambda)^{d-1}(\rho/\boldsymbol{r}\lambda + \partial F_t(c_t(0))/\partial a_t) = 0$ . Donc, si $\boldsymbol{c} = (c_t, \boldsymbol{c}_q), c_t$ est arbitraire et $\boldsymbol{c}_q = 0$.

Cherchons un vecteur propre à distance finie $\boldsymbol{c} + \boldsymbol{w}$ pour une v.p. $\boldsymbol{\varpi} = -\rho t/\boldsymbol{r}\lambda + u$. On va avoir :

$\boldsymbol{\varpi}(\boldsymbol{c} + \boldsymbol{w}) = (u - \rho t/\boldsymbol{r})(\boldsymbol{c} + \boldsymbol{w}) = \boldsymbol{x}\partial F(\boldsymbol{c}(0))/\partial a\, (\boldsymbol{c} + \boldsymbol{w})$

Donc : $u(\boldsymbol{c} + \boldsymbol{w}) - \rho t/\boldsymbol{r}\boldsymbol{w} = \boldsymbol{x}\partial F(\boldsymbol{c}(0))/\partial a\, \boldsymbol{w}$. Choisissons la composante $w_t$ telle que $w_t + c_t = 0$.

Projetons l'équation sur l'axe $t$ :

$-\rho/\boldsymbol{r}w_t = \partial F_t(\boldsymbol{c}(0))/\partial a_t\, w_t + \partial F_t(\boldsymbol{c}(0))/\partial a_q\, \boldsymbol{w}_q$. Comme $\rho/\boldsymbol{r} + \partial F_t(\boldsymbol{c}(0))/\partial a_t = 0$, on a :

$\partial F_t(\boldsymbol{c}(0))/\partial a_q\, \boldsymbol{w}_q = 0$. Projetons ensuite l'équation sur le sous-espace des $\boldsymbol{q}$ :

$\boldsymbol{\varpi}\boldsymbol{w}_q = \boldsymbol{q}\partial F_q(\boldsymbol{c}_t(0))/\partial a_t\, w_t + \boldsymbol{q}\partial F_q(\boldsymbol{c}_t(0))/\partial a_q\, \boldsymbol{w}_q$

Si l'on prend $w_t = c_t = 0$, on obtient une configuration spatiale : $\boldsymbol{\varpi}\boldsymbol{w}_q = \boldsymbol{q}\partial F_q(\boldsymbol{c}_t(0))/\partial a_q\, \boldsymbol{w}_q$ indépendante de $t$. $\boldsymbol{w}_q$ est vecteur propre de $\boldsymbol{q}\partial F_q/\partial a_q$ pour la v.p. $\boldsymbol{\varpi}$ avec $\partial F_t(\boldsymbol{c}_t(0))/\partial a_q\, \boldsymbol{w}_q = 0$.

**4- Remarques**

- Pour la commodité des calculs, on a pris 1 seule v.p. plus grande que 1, on aurait pu en prendre $q < d$.
- Si on a une seule variable, on a une EDO. Les résultats restent les mêmes.

- Il semble que la solution trouvée soit locale. En particulier, si l'on a plusieurs points fixes dans $C$, on a observé par des méthodes probabilistes, des passages d'un bassin de domination à l'autre.
- L'équation de Navier Stokes se met facilement sous la forme $\partial \boldsymbol{a} / \partial t = F(\boldsymbol{a})$ $[6]$. On a initialement $n+1$ inconnues $(u,p)$ avec autant de variables $(x,t)$ et d'équations différentielles. Notons $\partial u_i / \partial x = \boldsymbol{b}_i$, $i = 1,...,n$ $\Delta u_i = \Sigma_{j=1}^{n} \partial^2 u_i / \partial x_j^2 = \Sigma_{j=1}^{n} \partial \boldsymbol{b}_i / \partial x_j$, $\boldsymbol{c} = (x,t)$ $d_i = \partial p / \partial x_i$. On rajoute ainsi autant 'équations que de variables, :

$$\partial \boldsymbol{u} / \partial t + \Sigma_{j=1}^{n} u_j b_{ij} = \nu \Sigma_{j=1}^{n} \partial \boldsymbol{b}_i / \partial x_j - d_i + f_i(\boldsymbol{c})$$

$$\Sigma_{j=1}^{n} \partial \boldsymbol{u}_j / \partial x_j = 0 , \quad \partial u_i / \partial x = \boldsymbol{b}_i , \ i = 1,...,n , \ \boldsymbol{c} = (x,t)$$

Si les forces extérieures $f_i$ sont $C^\infty$ (sans doute, peut-on améliorer nos hypothèses) et si l'itération différentielle est bornée dans l'espace, alors la solution sera presque périodique comme le veut la proposition précédente. S'il est difficile de trouver des solutions périodiques, comme dans le problème $[8]$, en revanche, les fonctions presque périodiques résolvent la question. L'attracteur de Lorenz en est la preuve évidente. Les $f_i$ méritent d'être précisés dans l'énoncé du problème.
- On peut traiter les itérations avec décalage : $\boldsymbol{u}_1(\boldsymbol{x} + \boldsymbol{log}\lambda) = F(\boldsymbol{u}(\boldsymbol{x} + \boldsymbol{T}))$ : en posant $\boldsymbol{T} = \log(\tau)$, $\boldsymbol{x} = \log(\boldsymbol{t} / \tau)$ et $\boldsymbol{a}(\boldsymbol{t}) = \boldsymbol{u} \circ \log(\boldsymbol{t}))$, on a $\boldsymbol{a}_1(\lambda \boldsymbol{t} / \tau) = F(\boldsymbol{a}(\boldsymbol{t})$ et l'on remplace la v.p. $\lambda$ par $(\lambda\tau)$.
On peut traiter ainsi les EDO et EDP avec décalage sans plus de difficulté.

## VI. Références

## En guise de conclusion

Il s'agit ici d'une analyse rapide des itérations inversibles qui s'avère être un domaine délicat. Après vérification par des spécialistes, elle peut ouvrir aux théoriciens et aux praticiens de nouvelles perspectives. Car en ce domaine, des théories et des expérimentations, si rigoureuses soient-elles, ont souvent été conduites sans fil directeur capable de mener à une compréhension globale.
Les itérations avec de multiples inverses ont été étudiées précédemment dans un contexte probabiliste aboutissent à d'autres conclusions. Ces méthodes ont été appliquées aux itérations inversibles sans être sûr de la convergence vers leurs solutions. Le rapprochement de ces divers résultats reste essentiel.
D'ailleurs, dans l'espace des fonctions de Weierstrass où l'on travaille ici, l'écart est faible avec les processus stochastiques. Si cette théorie présentée dans un cadre simple est confirmée, il ne fait pas de doute que des généralisations vont être avancées par des spécialistes. Il est alors permis de conjecturer que même les courbes définies géométriquement sont représentables par des fonctions de Weierstrass.
Dans ce bref essai, on s'est attaché à démontrer les questions qui paraissaient essentielles, sans pouvoir les confronter à d'autres idées. Toutes les remarques et observations, pourvu qu'elles ne soient pas trop difficiles, seront donc les bienvenues.